\documentclass{amsart}
\usepackage[cp1250]{inputenc}
\usepackage{amsmath,amssymb,amsthm}
\newtheorem{definition}{Definition}[section]

\newtheorem{remark}{Remark}[section]

\newtheorem{theorem}{Theorem}[section]

\numberwithin{equation}{section}

\begin{document}
\title{Law Of Large Numbers For Random Dynamical Systems }
\author{Katarzyna Horbacz}
\address{Department of Mathematics, University of Silesia, Bankowa 14, 40-007 Katowice, Poland, (KH)}
\email{horbacz@math.us.edu.pl}
\author{ Maciej \'Sl\c eczka}
\address{Department of Mathematics, University of Silesia, Bankowa 14, 40-007 Katowice, Poland, (M\'S)}
\email{sleczka@math.us.edu.pl}

\begin{abstract}
We cosider random dynamical systems with randomly chosen jumps. The choice of deterministic dynamical system and jumps depends on a position. We proove the existence of an exponentially attractive invariant measure and the strong law of large numbers.
\end{abstract}

\maketitle

\section{Introduction}

In the present paper we are concerned with the problem of proving the law of large numbers (LLN) for random dynamical systems.

The question of establishing the LLN for an additive functional of a Markov process is one of the most fundamental in probability theory
and there exists a rich literature on the subject, see e.g. the monograph of
Meyn and Tweedie \cite{M1} and the citations therein. However, in most of the
existing results, 
it is usually assumed that the process under
consideration is stationary and its equilibrium state  is stable in some
sense, usually in the $L^2$, or total variation norm. Our stability condition is
formulated  in a weaker metric
than the total variation distance.

The law of large numbers we study in this note was also considered in many papers. Our results are based on a version of the  law of large numbers due to  Shirikyan (see \cite{S1}, \cite{S2}).
Recently Komorowski, Peszat and Szarek \cite{KPS} obtained the weak law of large numbers for the passive tracer model in a compressible environment and Walczuk studied Markov processes with the transfer operator having spectral gap in the Wasserstein  metric and proved the LLN in the non-stationary case \cite{W1}.

Random dynamical systems \cite{15}, \cite{HD}take into consideration some very important and widely
studied cases, namely dynamical systems generated by learning systems \cite{2}, \cite{22}, \cite{23},
\cite{35}, Poisson driven stochastic differential equations  \cite{17}, \cite{34}, \cite{48}, \cite{49},
iterated function systems with an infinite family of transformations \cite{30}, \cite{50}, \cite{51},
random evolutions \cite{12}, \cite{42} and
irreducible Markov systems \cite{52}.

A large class of applications of such models,
both in physics and biology, is worth mentioning here: the shot noise, the photo
conductive detectors, the growth of the size of structural populations, the motion
of relativistic particles, both fermions and bosons (see \cite{8},
\cite{24}, \cite{28}), the generalized stochastic
process introduced in the recent model of gene expression by Lipniacki
et al. \cite{Lip}.

A number of results have been obtained that claim an existence of an asymptotically stable, unique invariant measure for Markov processes generated by random dynamical systems for which the state space need not be locally compact. We consider  random dynamical systems with randomly chosen jumps acting on a given Polish space $(Y,\varrho)$.

The aim of this paper is to study stochastic processes whose paths follow deterministic dynamics between random times, jump times, at which they change their position randomly. Hence, we
analyse stochastic processes in
which  randomness appears at times  $t_0 < t_1 <t_2<\ldots $ We assume that a
point $ x_0 \in Y $ moves according to one of the dynamical systems $ \Pi_i :
{\mathbb{R}}_+ \times Y \to Y $ from some set $\{ \Pi_1, \ldots, \Pi_N \}$. 
The motion of the process is
governed by the equation $ X(t) = \Pi_i (t, x_0) $ until the first jump time
$t_1$. Then we choose a transformation $q_s : Y \to Y$  from  a family $\{q_s\, : \, s\in S = \{1, \ldots , K\}\} $ and define $x_1 = q_s(\Pi_i (t_1, x_0))$.  The process restarts from that
new point $x_1$ and continues as before. This gives the stochastic process
$\{X(t)\}_{t \ge 0}$ with jump times  $\{t_1, t_2, \ldots \}$ and post jump
positions $\{x_1, x_2, \ldots \}$. The probability determining the frequency with
which the dynamical systems $\Pi_i$  are chosen is described by a  matrix of
probabilities ${[p_{ij}]}_{i,j=1}^N $, $p_{ij} : Y \to [0, 1]$. The maps  $q_s$ are randomly chosen with place dependent distribution. Given a Lipschitz function $\psi:X \to \mathbb{R}$ we define
$$ 
S_n(\psi)= \psi(x_0)+\dots+\psi(x_n).
$$
Our aim is to find conditions under which $S_n(\psi)$ satisfies law of large numbers.
Our results are based on an exponential convergence theorem due 
to  \'Sl\c eczka and  Kapica (see \cite{KS}) and a version of the law of large numbers due to  Shirikyan (see \cite{S1}, \cite{S2}).

\section{Notation and basic definitions}

Let $(X,d)$ be a \emph {Polish} space, i.e. a complete and separable metric space and denote by $\mathcal{B}_X$ the 
$\sigma$-algebra of Borel subsets of $X$. By $B_b(X)$ we denote the space of bounded Borel-measurable functions equipped with the supremum norm, $C_b(X)$ stands for the subspace of bounded continuous functions. Let 
$\mathcal{M}_{fin}(X)$ and $\mathcal{M}_1(X)$ be the sets of Borel measures on $X$ such that $\mu(X)<\infty$ for $\mu\in\mathcal{M}_{fin}(X)$ and $\mu(X)=1$ for $\mu\in\mathcal{M}_1(X)$. The elements of $\mathcal{M}_1(X)$ are called \emph {probability measures}. The elements of $\mathcal{M}_{fin}(X)$ for which  $\mu(X)\le 1$ are called \emph {subprobability measures}. By $supp\, \mu$ we denote the support of the measure $\mu$. We also define
$$
{\mathcal M}_1^L (X)=\{\mu\in{\mathcal M}_1(X):\,\int_X L(x)\mu (d x)<\infty\}
$$
where $L:X\to [0,\infty)$ is an arbitrary Borel measurable function and
$$
{\mathcal M}_1^1 (X)=\{\mu\in{\mathcal M}_1(X):\,\int_X d(\bar x ,x)\mu (d x)<\infty\},
$$
where $\bar x\in X$ is fixed. By the triangle inequality this family is independent of the choice of $\bar x$.
\newline
The space $\mathcal{M}_1(X)$ is equipped with the \emph{Fourtet-Mourier metric}:
$$
\Vert\mu_1-\mu_2\Vert_{FM}=\sup\{|\int_X f(x)(\mu_1-\mu_2)(dx)|:\, f\in\mathcal{F}\},
$$
where
$$
\mathcal{F}=\{f\in C_b(X):\, |f(x)-f(y)|\le d(x,y) \quad\text{and}\quad |f(x)|\le 1\quad\text{for}\quad x,y\in X\}.
$$
Let $P:B_b(X)\to B_b(X)$ be a \emph{Markov operator}, i.e. a linear operator satisfying $P{\bf 1}_X={\bf 1}_X$ and $Pf(x)\ge 0$ if $f\ge 0$. Denote by $P^{*}$ the the dual operator, i.e operator $P^{*}:\mathcal{M}_{fin}(X)\to
\mathcal{M}_{fin}(X)$ defined as follows 
$$
P^{*}\mu(A):=\int_X P {\bf 1}_A(x)\mu(dx)\qquad\text{for}\qquad A\in\mathcal{B}_X.
$$
We say that a measure $\mu_*\in\mathcal{M}_1(X)$ is \emph{invariant} for $P$ if
$$
\int_X Pf(x)\mu_*(dx)=\int_X f(x)\mu_*(dx)\qquad\text{for every}\qquad f\in B_b(X)
$$
or, alternatively, we have $P^* \mu_*=\mu_*$. An invariant measure $\mu$ is \emph{attractive} if 
$$
\lim\limits_{n\to\infty}\int_X P^nf(x)\, \mu(dx)=\int_X f(x)\, \mu(dx)\quad\text{for}\quad f\in C_b(X),\, \mu\in\mathcal{M}_1(X).
$$
\newline
By $\{\mathbf{P}_x:\, x\in X\}$ we denote a \emph{transition probability function} for $P$, i.e. a family of measures $\mathbf{P}_x\in\mathcal{M}_1(X)$ for $x\in X$, such that the map $x\mapsto\mathbf{P}_x(A)$ is measurable for every $A\in\mathcal{B}_X$ and 
$$
Pf(x)=\int_X f(y) \mathbf{P}_x(dy)\qquad\text{for}\qquad x\in X\quad\text{and}\quad f\in B_b(X)
$$
or equivalently $P^*\mu(A)=\int_X \mathbf{P}_x(A)\mu(dx)$ for $A\in\mathcal{B}_X$ and $\mu\in\mathcal{M}_{fin}(X)$. 
We say that a  vector $(p_1, \ldots ,p_N)$ where $p_i :Y \to [0, 1]$ is {\it a probability vector} if
$$
\sum_{i=1}^N p_i(x) = 1 \quad \textrm{for}\quad x\in Y.
$$
Analogously a matrix $[p_{ij}]_{i,j }$ where $ p_{ij} :Y \to [0, 1]$ for $i, j \in \{1, \ldots, N\}$
is {\it a probability matrix } if
$$
\sum_{j=1}^N p_{ij}(x) = 1  \quad \textrm{for}\quad x\in Y\quad \textrm{and} \quad i \in \{1, \ldots, N \}.
$$

\begin{definition}
A coupling for $\{\mathbf{P}_x: x\in X\}$ is a family $\{\mathbf{B}_{x,y}:\, x,y\in X\}$ of probability measures on $X\times X$ such that for every $B\in\mathcal{B}_{X^2}$ the map $X^2\ni (x,y)\mapsto \mathbf{B}_{x,y}(B)$ is measurable and
$$ 
\mathbf{B}_{x,y}(A\times X)=\mathbf{P}_x(A),\qquad \mathbf{B}_{x,y}(X\times A)=\mathbf{P}_y(A)
$$
for every $x,y\in X$ and $A\in\mathcal{B}_X$.
\end{definition}
In the following we assume that there exists a \emph{subcoupling} for $\{\mathbf{P}_x: x\in X\}$, i.e. a family $\{\mathbf{Q}_{x,y}:\,x,y\in X\}$ of subprobability measures on $X^2$ such that the map $(x,y)\mapsto \mathbf{Q}_{x,y}(B)$ is measurable for every Borel $B\subset X^2$ and
$$
\mathbf{Q}_{x,y}(A\times X)\le \mathbf{P}_x(A)\qquad\text{and}\qquad \mathbf{Q}_{x,y}(X\times A)\le \mathbf{P}_y(A)
$$
for every $x,y\in X$ and Borel $A\subset X$.
\newline
Measures $\{\mathbf{Q}_{x,y}:x,y\in X\}$ allow us to construct a coupling for $\{\mathbf{P}_x:x\in X\}$. Define on $X^2$ the family of measures $\{\mathbf{R}_{x,y}:x,y\in X\}$ which on rectangles $A\times B$ are given by
$$
\mathbf{R}_{x,y}(A\times B)=\frac{1}{1-\mathbf{Q}_{x,y}(X^2)}(\mathbf{P}_x(A)-\mathbf{Q}_{x,y}(A\times X)) (\mathbf{P}_y(B)-\mathbf{Q}_{x,y}(X\times B)),
$$
when $\mathbf{Q}_{x,y}(X^2)<1$ and $\mathbf{R}_{x,y}(A\times B)=0$ otherwise. A simple computation shows that the family $\{\mathbf{B}_{x,y}:\, x,y\in X\}$ of measures on $X^2$ defined by
\begin{equation}\label{BQR}
\mathbf{B}_{x,y}=\mathbf{Q}_{x,y}+\mathbf{R}_{x,y}\quad\text{for}\quad x,y\in X 
\end{equation}
is a coupling for $\{\mathbf{P}_x:\, x\in X\}$.
\newline
The following Theorem due to M. \'Sl\c eczka and R. Kapica (see \cite{KS}) will be used in the proof of Theorem \ref{HS1} in Section 4.
\begin{theorem}\label{KS} Assume that a Markov operator $P$  and transition probabilities $\{\mathbf{Q}_{x,y}:\,x,y\in X\}$ satisfy
\newline
{\bf A0} $P$ {\it is a} \emph{ Feller operator}, {\it i.e.} $P(C_b(X))\subset C_b (X)$. 
\newline
{\bf A1} {\it There exists a} \emph{ Lapunov function}\, {\it for} $P$, {\it i.e. continuous function} $L: X\to [0, \infty )$     {\it such that} 
$L$ {\it is bounded on bounded sets, }
$\lim _{x\to \infty}L(x)=+\infty$ {\it and for some} $\lambda \in (0,1), \, c>0$
$$
PL (x)\le \lambda L(x) +c \qquad for \qquad x\in X.
$$
{\bf A2} {\it There exist} $F\subset X^2$ {\it and} $\alpha \in  (0,1)$ {\it such that} $supp \,{\mathbf Q}_{x,y}\subset F$ {\it and}
\begin{equation}\label{A2}
\int_{ X^2} d (u,v){\mathbf Q}_{x,y}(d u,d v) \le \alpha d(x,y)\qquad for \qquad (x,y)\in F.
\end{equation}
{\bf A3} {\it There exist} $\delta>0,\, l>0$ {\it and} $\nu\in (0,1]$ {\it such that} 
$$
1- {\mathbf Q}_{x,y}(X^2)  \le l d(x,y)^{\nu}
$$
and
$$
{\mathbf Q}_{x,y} (\{(u,v)\in X^2:\, d(u,v)<\alpha d(x,y)\} )\ge\delta
$$
for $(x,y)\in F$
\newline
{\bf A4} {\it There exist} $\beta\in (0,1)$, ${\tilde C}>0$ {\it and} $R>0$ {\it such that for} 
$$
\kappa (\,(x_n,y_n)_{n\in{\mathbb N}_0}\,)=\inf \{n\in {\mathbb N}_0 :\, (x_n,y_n)\in F\quad\text{and}\quad L(x_n)+L(y_n)<R\}
$$
{\it we have}
$$
{\mathbb E}_{x,y} \beta^{-\kappa}\le {\tilde C}\qquad whenever \qquad L(x)+L(y)<\frac{4c}{1-\lambda},
$$
{\it where} $\mathbb{E}_{x,y}$ {\it denotes here the expectation with respect to the chain starting from} $(x,y)$ {\it and with trasition function } $\{\mathbf{B}_{x,y}:\, x,y\in X\}$. \newline
Then operator $P$ possesses a unique invariant measure $\mu_{*}\in\mathcal{M}_1^L (X)$, which is attractive in $\mathcal{M}_1(X)$.
Moreover, there exist $q\in (0,1)$ and $C>0$ such that
\begin{equation}\label{mix}
\| P^{* n}\mu -\mu_{*}\|_{FM}\le q^n C(1+\int_X L(x)\mu (dx) )
\end{equation}
for $\mu\in\mathcal{M}_1^L (X)$ and $n\in\mathbb{N}$.
\end{theorem}
We will also need a version of the strong law of large numbers due to A. Shirikyan (\cite{S1}, \cite{S2}). It is originally 
formulated for Markov chains on a Hilbert space, however analysis of the proof shows that it remains true for Polish spaces.

\begin{theorem}\label{Shir} Let $(\Omega,\mathcal{F},\mathbb{P})$ be a probability space and let $X$ be a Polish space. 
Suppose that for a family of Markov chains $((X_n^x)_{n\ge 0},\mathbb{P}_x)_{x\in X}$ on $X$ with Markov operator 
$P:B_b(X)\to B_b(X)$ there exists a unique invariant measure $\mu_*\in\mathcal{M}_1(X)$, a continuous function $v:X\to\mathbb{R}_+$ and a sequence $(\gamma_n)_{n\in\mathbb{N}}$ of positive numbers such that 
$\gamma_n\to 0$ as $n\to\infty$ and
$$
||P^{*n}\delta_x-\mu_*||_{FM}\le \gamma_n v(x)\quad\text{for}\quad x\in X.
$$
If
$$
C=\sum_{n=0}^{\infty}\gamma_n <\infty
$$
and there exits a continuous function $h:X\to\mathbb{R}_+$ such that 
$$
\mathbb{E}_x(v(X_n^x))\le  h(x)\quad\text{for}\quad x\in X,n\ge 0,
$$
where $\mathbb{E}_x$ is the expectation with respect to $\mathbb{P}_x$, then for any $x\in X$ and any 
bounded Lipschitz function $f:X\to \mathbb{R}$ we have
$$
\lim _{n\to\infty}\frac{1}{n}\sum_{k=0}^{n-1} f(X_k^x)=\int_X f(y)\, \mu_*(dy)
$$
$\mathbb{P}_x$ almost surely.
\end{theorem}

\section{ Random Dynamical Systems}
Let  $(Y, \varrho)$ be a Polish space,  $\mathbb{R}_+=[0,+\infty )$ and $I = \{1, \dots ,N\}$, $ S = \{1, \ldots , K\}$,
where $N$ and $K$ are given positive integers.

Let  $\Pi_{i}:\mathbb{R}_+\times Y\rightarrow Y$, $i \in I$,  be a finite sequence of semidynamical  systems, i.e.
$$
\Pi_i(0,x)=x \quad \text{for }\quad i\in I, \,\, x\in Y
$$
and
$$
\Pi_i(s+t,x)=\Pi_i(s,(\Pi_i(t,x)) \quad \textrm{for }\quad s,t \in \mathbb{R}_+, \,\,i\in I \,\,\,\textrm{and}\,\,\, x\in Y.
$$

We are given probability vectors $p_i :Y \to [0,1], \,\,\, i \in I$, $\overline {p}_s :Y \to [0,1],
\,\,\, s \in S $, a matrix of probabilities  $[p_{ij}]_{i, j \in I}$, $ p_{ij}:Y\rightarrow [0, 1], \,\,\, i,j \in I$
and a family of continuous functions $q_s : Y \to Y , s \in S $. In the sequel we denote the system by $(\Pi , q, p )$.

Finally, let  $(\Omega, \Sigma, \mathbb{P} )$ be a probability space  and   $\{t_n\}_{n\ge 0}$ be an increasing
 sequence of random variables
$t_n :\Omega \to \mathbb{R}_+$ with $t_0 =0$   and such that the increments
$\Delta t_n=t_n-t_{n-1}$,   $n \in \mathbb{N} $,
 are independent  and have the same density  $g(t)=\lambda e^{-\lambda t}$, $ t \ge 0 $.

The action of randomly chosen dynamical systems, with randomly chosen  jumps, at random moments $t_k$ corresponding to the system $ (\Pi , q, p )$ can be roughly described as follows.

We choose an initial point $x_0 \in Y $ and randomly select a transformation $\Pi_i$ from the set $\{\Pi_1 , \ldots , \Pi_N \}$ in such a way that the probability of choosing $\Pi_i$ is equal to $p_i(x_0)$, and we define
$$
X(t) = \Pi_i(t, x_0) \quad \textrm{for}\quad 0\le t < t_1.
$$
Next, at the random moment $t_1$, at the point $\Pi_i(t_1, x_0)$ we choose a jump $q_s$ from the set $\{q_1, \ldots ,q_K\}$ with probability $\overline {p}_s(\Pi_i(t_1, x_0 ))$. Then we define
$$
x_1 = q_s (\Pi_i (t_1, x_0)).
$$
After that we choose $\Pi_{i_1}$ with probability $p_{ii_1}(x_1)$,  define
$$
X(t) = \Pi_{i_1} (t - t_1, x_1 )\quad \textrm{for}\quad t_1 < t <t_2
$$
and at the point $\Pi_{i_1}(t_2 - t_1, x_1 )$ we choose $q_{s_1}$ with probability $\overline {p}_{s_1}(\Pi_{i_1}(t_2 - t_1, x_1))$. Then we define
$$
x_2 = q_{s_1}(\Pi_{i_1} (t_2 - t_1, x_1 )).
$$

Finally, given $x_n$, $n\ge 2 $, we choose $ \Pi_{i_n} $ in such a way that the probability of  choosing  $ \Pi_{i_n} $ is equal to $p_{i_{n-1}i_n}(x_n)$ and we define
$$
X(t) = \Pi_{i_n} (t - t_n, x_n )\quad \textrm{for}\quad t_n < t <t_{n+1}.
$$
At the point  $ \Pi_{i_n}(\Delta t_{n+1},  x_n ) $ we choose $q_{s_n}$ with probability  $\overline {p}_{s_n}(\Pi_{i_n}(\Delta t_{n+1}, x_n))$. Then we define
$$
x_{n+1} = q_{s_n}(\Pi_{i_n} (\Delta t_{n+1}, x_n )).
$$

We obtain a piecewise-deterministic trajectory for $\{X(t)\}_{t \ge 0}$ with jump
times  $\{t_1, t_2, \ldots \}$ and post jump locations $\{x_1, x_2, \ldots \}$.

We may reformulate  the above considerations  as follows: Let $\{\xi_n\}_{n \ge 0}$ and
$\{\eta _n\}_{n \ge  1}$ be sequences of random variables, $\xi_n :\Omega \to I$ and
$ \eta_n :\Omega \to S $   and let $\{y_n\}_{n \ge 1}$ be auxiliary random  variables,
$y_n : \Omega \to Y $,  such that
\begin{equation}\label{(4.1.1)}
\begin{aligned}
&\mathbb{P} (\xi_0 = i | x_0 = x ) = p_i (x),\\
&\mathbb{P} (\xi_n = k | x_n = x \quad \textrm{and} \quad \xi_{n-1} = i ) = p_{ik}(x),\\
\end{aligned}
\end{equation}
and
\begin{equation}\label{(4.1.2)}
\begin{aligned}
& y_n = \Pi_{\xi_{n-1}} (t_n - t_{n-1}, x_{n-1}),\\
&\mathbb{P} (\eta_n = s | y_n = y ) = \overline{p}_s (y)
\end{aligned}
\end{equation}
for $ n\ge 1, \,\,x, y \in Y,\,\, k, i \in I $ and $ s \in S $ .

Assume that $\{\xi_n\}_{n \ge 0}$ and $\{\eta_n\}_{n \ge 0} $ are independent  of $\{t_n\}_{n \ge 0}$ and that for every $n \in \mathbb{N}$ the variables $\eta_1, \ldots ,\eta_{n-1}$, $ \xi_1, \ldots ,\xi_{n-1}$ are also independent.

Given an initial random variable $\xi_0$ the sequence of the random variables  $\{x_n\}_{n\ge 0}$,  $x_n : \Omega \to Y $, is given by  
\begin{equation}\label{(4.1.3)}
x_n =q_{\eta_n} \big (\Pi_{\xi_{n-1}}(t_n - t_{n-1}, x_{n-1})\big ) \quad
\text{for}\quad n=1,2, \dots
\end{equation}
and the stochastic process $\{X(t)\}_{t \ge 0}$, $X(t) : \Omega \to Y$, is given by
\begin{equation}\label{(4.1.4)}
X(t) = \Pi_{\xi_{n-1}}(t - t_{n-1}, x_{n-1}) \quad \textrm{for} \quad t_{n-1} \le
t < t_n,\quad n = 1,2, \ldots
\end{equation}

It is easy to see that $\{X(t)\}_{t \ge 0}$ and $\{x_n\}_{n \ge 0}$ are not Markov processes. In order to
use the  theory of Markov operators we must  redefine the processes $\{X(t)\}_{t \ge 0}$ and $\{x_n\}_{n \ge 0}$ in such a way   that the redefined processes become Markov.

For this purpose, consider the space $Y\times I $ endowed with the metric $d$ given by
\begin{equation}\label{(4.1.5)}
d \big((x, i), (y, j)\big)=\varrho(x,y) + \varrho_d(i,
j)\quad\textrm{for}\quad x, y\in Y, \,\,i, j\in I,
\end{equation}
where $\varrho_d$ is the discrete metric in $I$.\newline
Now  define a stochastic process $\{\xi (t) \}_{t \ge 0}$, $\xi (t): \Omega \to I $, by
$$
\xi (t) = \xi_{n-1} \quad \textrm{for} \quad  t_{n-1} \le t <t_{n},\quad n=1, 2,  \ldots
$$
Then the stochastic process $ \{(X(t), \xi (t))\}_{t \ge 0}$,  $(X(t), \xi (t)) : \Omega \to Y \times I $ has the required Markov property.

In many applications we are mostly interested in values of the process $X(t)$  at the switching
points $t_n $. Therefore, we will also study the stochastic discrete process (post jump locations)
$\{(x_n, \xi_n) \}_{n\ge 0}$ , $(x_n, \xi_n) : \Omega \to Y \times I $. Clearly
$\{(x_n, \xi_n) \}_{n \ge 0}$  is a Markov process too.

We consider the stochastic process  $\{(x_n, \xi_n) \}_{n\ge 0}$ , $(x_n, \xi_n) : \Omega \to Y
\times I $, defined by \eqref{(4.1.1)}--\eqref{(4.1.3)} with the help of the system $(\Pi, q, p)$.
We will need the following assumptions:

The transformations  $\Pi_i : \mathbb{R}_+ \times Y \to Y$, $i\in I$ and $q_s :
Y \to Y $, $s\in S $, are continuous and  there exists $x_* \in Y$ such that
\begin{equation}\label{(4.2.1)}
\int_{\mathbb{R}_+}e^{-\lambda t} \varrho (q_s (\Pi_j(t,x_*)) , q_s(x_*))\ dt < \infty
\quad \textrm{for} \quad j \in I, \quad s \in S.
\end{equation}
The functions $\overline{p}_s$, $s \in S$, and $p_{ij}$,  $i,j \in I$,   satisfy the following conditions
\begin{equation}\label{(4.2.2)}
\begin{split}
\sum_{j\in I} |p_{ij}(x) - p_{ij}(y)| &\le L_p\varrho (x,y)  \quad
\textrm{for}\quad x,y \in Y, \,\, i \in I,\\
\sum_{s\in S} |\overline{p}_{s}(x) - \overline{p}_{s}(y)| &\le L_{\overline{p}}\varrho (x,y)
\quad \textrm{for}\quad x,y \in Y,
\end{split}
\end{equation}
where  $L_p, L_{\overline{p}} > 0$. 

We also assume  that for the system $(\Pi , q, p)$ there are three constants $L\ge 1$, $\alpha \in \mathbb{R} $ and $L_q > 0$ such that
\begin{equation}\label{(4.2.4)}
\sum_{j\in I} p_{ij}(y)\varrho (\Pi_j(t,x) ,\Pi_j(t,y)) \le Le^{ \alpha t}\varrho (x,y) \quad
\textrm{for}\quad x,y \in Y, \,\, i \in I,  \,\,  t \ge 0
\end{equation}
and
\begin{equation}\label{(4.2.5)}
\sum_{s \in S} \overline{p}_s(x)\varrho (q_s(x),q_s(y))  \le L_q \varrho (x,y)  \quad
\textrm{for} \quad x,y \in Y.
\end{equation}

For $x, y \in Y,\, t\ge 0$ we define
\begin{equation}
\begin{aligned}
&I_{\Pi}(t, x, y) = \{ j \in I\quad: \varrho (\Pi_j(t, x) , \Pi_j(t, y)) \le Le^{\alpha t}\varrho (x , y)\}\\
&I_{q}(x, y) = \{ s \in S \quad: \varrho (q_s(x) , q_s(y)) \le L_q\varrho (x , y)\}
\end{aligned}
\end{equation}

Assume that there are $p_0 > 0, q_0 > 0$ such that : for every $i_1, i_2 \in I, \, x, y \in Y$ and $ t \ge 0$ we have

\begin{equation}\label{(S1)}
\begin{aligned}
&\sum_{j \in I_{\Pi}(t, x, y)} p_{i_1j}(x)p_{i_2j}(y) > p_0,\\
&\sum_{s \in I_{q}(x, y)}\overline{p}_{s}(x)\overline{p}_{s}(y) > q_0.
\end{aligned}
\end{equation}

\begin{remark}
The condition \eqref{(S1)} is satisfied if there are $i_0 \in I, s_0 \in S$ such that 
\begin{equation}
\begin{aligned}
&\varrho (\Pi_{i_0}(t,x) , \Pi_{i_0}(t,y)) \le Le^{ \alpha t}\varrho(x,y) \quad
\textrm{for}\quad x,y \in Y, \,\,   t \ge 0,\\
&\varrho (q_{s_0}(x),q_{s_0}(y))  \le L_q \varrho (x,y)  \quad
\textrm{for} \quad x,y \in Y,
\end{aligned}
\end{equation}
and
\begin{equation}
\begin{aligned}
&\inf_{i\in I} \inf_{x \in Y} p_{i i_0}(x) > 0,\\
&\inf_{x \in Y} \overline{p}_{s_0}(x) > 0.
\end{aligned}
\end{equation}

\end{remark}

To begin our study of the stochastic process  $\{(x_n, \xi_n) \}_{n\ge 0}$ consider the  sequence of distributions
$$
\overline\mu_n(A) = \mathbb{P} \big ((x_n, \xi_n) \in A \big )  \quad \textrm{for} \quad A \in \mathcal{B} (Y\times I), \, n \ge 0.
$$
It is easy to see that there exists  a Markov-Feller operator $P: \mathcal{M} \to \mathcal{M} $ such that
\begin{equation*}\label{(4.2.6)}
\overline\mu_{n+1} = P \overline\mu_n \quad \textrm{for } \quad n \ge 0.
\end{equation*}
The operator $P$ is given by the formula
\begin{equation}\label{operatorP}
P\mu (A) = \sum_{j\in I} \sum_{s \in S} \int_{Y\times I}   \int_0^{+\infty} \lambda e^{-\lambda t}
 1_A\big (q_s\big ( \Pi_j(t,x)\big ),j \big ) p_{ij}(x)\overline{p}_s\big(\Pi_j (t, x)\big) \, dt\,
\mu (dx, di)
\end{equation}
and its dual operator $U$  by 
\begin{equation}\label{operatorU}
Uf(x, i) =  \sum_{j\in I} \sum_{s \in S} \int_0^{+\infty }   \lambda e^{-\lambda
t} f\big (q_s\big ( \Pi_j(t, x )\big ),j\big )p_{ij}(x)\overline{p}_s\big(\Pi_j
(t, x)\big)  \,dt,
\end{equation}
where $\lambda$ is the intensity of
 the Poisson process which governs the increment $\Delta t_n$ of the random variables $\{t_n\}_{n\ge 0}$.
The operator $P$ given by \eqref{operatorP} is called a {\it  transition operator}
for this system.

\section{The main theorem}

\begin{theorem}\label{HS1}
Assume that system $(\Pi, p, q)$ satisfies conditions \eqref{(4.2.1)}--\eqref{(S1)}. If
\begin{equation}\label{(S2)}
LL_q + \frac{\alpha}{\lambda} < 1.
\end{equation}
then \newline
(i) there exists a unique invariant measure $\mu_*\in\mathcal{M}_1^1(Y\times I)$ for the process $(x_n,\xi_n)_{n\ge 0}$, 
which is attractive in $\mathcal{M}_1(Y\times I)$.\newline
(ii) there exist $q\in (0,1)$ and $C>0$ such that for $\mu\in\mathcal{M}_1^1(Y\times I)$ and $n\in\mathbb{N}$
$$
||P^{*n}\mu-\mu_* ||_{FM}\le q^n C(1+\int_Y \varrho (x,x_*) \,\mu(dx)),
$$
where $x_*$ is given by (\ref{(4.2.1)}),\newline
(iii) the strong law of large numbers holds for the process $(x_n,\xi_n)_{n\ge 0}$ starting from $(x_0,\xi_0 )\in Y\times I$, i.e. for every bounded Lipschitz function 
$f:Y\times I\to \mathbb{R}$ and every $x_0\in Y$ and $\xi_0\in I$ we have
$$
\lim_{n\in \infty}\frac{1}{n}\sum_{k=0}^{n-1} f(x_k,\xi_k)=\int_{Y\times I} f(x,\xi)\, \mu_*(dx,d\xi)
$$
$\mathbb{P}_{x_0,\xi_0}$ almost surely.
\end{theorem}

{\it Proof of Theorem \ref{HS1}}\newline
We are going to verify assumptions of Theorem \ref{KS}. Set $X=Y\times I$, $F=X\times X$ and define
$$
\begin{aligned}
&{\mathbf Q}_{(x_1,i_1)(x_2,i_2)}(A)=
\\ &\sum_{j\in I}\sum_{s\in S}\int_0^{+\infty}\lambda e^{-\lambda t}\{ p_{i_1 j}(x_1)
\overline{p}_s\big ( \Pi_j(t,x_1)\big )\wedge p_{i_2 j}(x_2)\overline{p}_s\big (\Pi_j(t,x_2)\big)\} \times \\
& \times 1_A\big (\big( q_s\big(\Pi_j(t,x_1)\big ),j),(q_s\big (\Pi_j(t,x_2)\big),j\big )\big )\, dt
\end{aligned}
$$
for $A\subset X\times X$, where $a\wedge b$ stands for the minimum of $a$ and $b$.\newline

{\bf A0}. The continuity of functions $p_{ij}, \overline{p}_s, q_s$ implies that the operator $P$ defined in 
(\ref{operatorP}) is a Feller operator.\newline

{\bf A1}. Define $L(x,i)=\varrho(x,x_*)$ for  $(x,i)\in X$. By (\ref{operatorU}) we have
$$
\begin{aligned}
UL(x,i)\le &\sum_{j\in I}\sum_{s\in S} \int_0^{+\infty} \varrho (q_s\big (\Pi_j(t,x)\big ) , q_s\big (\Pi_j(t,x_*)\big ))
\lambda e^{-\lambda t} p_{ij}(x)\overline{p}_s\big (\Pi_j(t,x)\big )\, dt \\
& +\sum_{j\in I}\sum_{s\in S} \int_0^{+\infty} \varrho (q_s\big (\Pi_j(t,x_*)\big ) , q_s(x_*))
\lambda e^{-\lambda t} p_{ij}(x)\overline{p}_s\big (\Pi_j(t,x)\big )\, dt \\
& +\sum_{j\in I}\sum_{s\in S} \int_0^{+\infty} \varrho (q_s(x_*),x_*)
\lambda e^{-\lambda t} p_{ij}(x)\overline{p}_s\big (\Pi_j(t,x)\big )\, dt .
\end{aligned}
$$
Further, using (\ref{(4.2.1)}), (\ref{(4.2.4)}) and (\ref{(4.2.5)}) we obtain 
\begin{equation}\label{Lap}
UL(x,i)\le aL(x,i)+b,
\end{equation}
where
\begin{equation}\label{constants}
\begin{aligned}
&a=\frac{\lambda L L_q}{\lambda-\alpha}, \\
&b=\sum_{j\in I}\sum_{s\in S}\int_0^{+\infty}\lambda e^{-\lambda t}\varrho (q_s\big (\Pi_j(t,x_*)\big ),q_s(x_*))\, dt 
+\sum_{s\in S}\varrho(q_s(x_*),x_*), 
\end{aligned}
\end{equation}

so $L$ is a Lapunov function for $P$.\newline

{\bf A2}. Observe that by (\ref{(4.1.5)}), (\ref {(4.2.4)}) and (\ref{(4.2.5)}) we have for $(x_1,i_1),(x_2,i_2)\in X$
$$
\begin{aligned}
&\int_{X^2} d(u,v)\,{\mathbf Q}_{(x_1,i_1)(x_2,i_2)}(du,dv)=\\
&\sum_{j\in I}\sum_{s\in S} \int_0^{+\infty}\lambda e^{-\lambda t}\{ p_{i_1 j}(x_1)
\overline{p}_s\big ( \Pi_j(t,x_1)\big )\wedge p_{i_2 j}(x_2)\overline{p}_s\big (\Pi_j(t,x_2)\big)\}\times \\
&\times \varrho (q_s\big (\Pi_j(t,x_1)\big),q_s\big (\Pi_j (t,x_2)\big ))\,dt \\
&\le \sum_{j\in I}\sum_{s\in S} \int_0^{+\infty}\lambda e^{-\lambda t} p_{i_1 j}(x_1)
\overline{p}_s\big ( \Pi_j(t,x_1)\big ) \varrho (q_s\big (\Pi_j(t,x_1)\big),q_s\big (\Pi_j (t,x_2)\big ))\,dt \\
& \le \beta \varrho (x_1,x_2) \le \beta \,d\big((x_1,i_1),(x_2,i_2)\big )
\end{aligned}
$$
with $\beta=\frac{\lambda L L_q}{\lambda-\alpha}<1$ by (\ref{(S2)}).\newline

{\bf A3}. From (\ref{(4.2.2)}) and (\ref{(4.2.4)}) it follows that
$$
\begin{aligned}
&1-\sum_{j\in I}\sum_{s\in S} \{ p_{i_1 j}(x_1)
\overline{p}_s\big ( \Pi_j(t,x_1)\big )\wedge p_{i_2 j}(x_2)\overline{p}_s\big (\Pi_j(t,x_2)\big)\} \\
&\le \sum_{j\in I}\sum_{s\in S} | p_{i_1 j}(x_1)
\overline{p}_s\big ( \Pi_j(t,x_1)\big ) - p_{i_2 j}(x_2)\overline{p}_s\big (\Pi_j(t,x_2)\big)| \\
& \le \sum_{j\in I}\sum_{s\in S} p_{i_1 j}(x_1)|\overline{p}_s\big (\Pi_j(t,x_1)\big )-\overline{p}_s\big (\Pi_j(t,x_2)\big )| \\
& + \sum_{j\in I}\sum_{s\in S} \overline{p}_s\big (\Pi_j(t,x_2)\big )|p_{i_1 j}(x_1)-p_{i_2 j}(x_2)| \\
&\le L L_{\overline {p}} e^{\alpha t}\varrho (x_1,x_2) + L_p \varrho (x_1,x_2) +2N \varrho_d(i_1,i_2)
\end{aligned}
$$
and consequently
$$
1-{\mathbf Q}_{(x_1,i_1)(x_2,i_2)}(X^2)\le (L_p+\frac{\lambda L L_{\overline{p}}}{\lambda-\alpha})\varrho (x_1,x_2)+2N\varrho_d(i_1,i_2).
$$
Fix $x_1,x_2\in Y$ and $i_1,i_2\in I$. Define $B=\{\big ( (u,j),(v,j)\big ):\, \varrho (u,v)<\beta \varrho (x_1,x_2), j\in I\}$. 
If $\alpha\ge 0$ then there exists $T_0>0$ such that $L L_q e^{\alpha t}<\beta$ for $t<T_0$. Set $A=(0,T_0)$. 
If  $\alpha <0$ then there exists $T_0>0$ such that $L L_q e^{\alpha t}<\beta$ for $t>T_0$. Set $A=(T_0,\infty)$. 
In both cases define $r=\int_A \lambda e^{-\lambda t}\, dt$. For all $x,y\in Y$, $t\in A$, $j\in I_{\Pi}(t,x,y)$ and 
$s\in I_q\big (\Pi_j(t,x),\Pi_j(t,y)\big )$ we have 
\begin{equation}\label{S_set}
\big ( (q_s(\Pi_j(t,x) ),j ), (q_s (\Pi_j(t,y) ),j )\big )\in B.
\end{equation}
From  (\ref{(S1)}) and (\ref{S_set}) we obtain 
$$
\begin{aligned}
&{\mathbf Q}_{(x_1,i_1)(x_2,i_2)} (B) \\
&\ge \int_A \lambda e^{-\lambda t}\sum_{j\in I_{\Pi}(t,x_1,x_2)}\sum_{s\in I_q (\Pi_j(t,x),\Pi_j(t,y) )} 
\{ p_{i_1 j}(x_1)\overline{p}_s\big ( \Pi_j(t,x_1)\big )\wedge p_{i_2 j}(x_2)\overline{p}_s\big (\Pi_j(t,x_2)\big)\} \times \\
&\times 1_B\big (\big ( q_s\big (\Pi_j(t,x_1)\big ),j\big ) ,\big(q_s\big (\Pi_j(t,x_2),j\big )\big ) \big)\, dt \\
&= \int_A \lambda e^{-\lambda t}\sum_{j\in I_{\Pi}(t,x_1,x_2)}\sum_{s\in I_q (\Pi_j(t,x),\Pi_j(t,y) )} 
\{ p_{i_1 j}(x_1)\overline{p}_s\big ( \Pi_j(t,x_1)\big )\wedge p_{i_2 j}(x_2)\overline{p}_s\big (\Pi_j(t,x_2)\big)\} \, dt\\
&\ge \int_A \lambda e^{-\lambda t}\sum_{j\in I_{\Pi}(t,x_1,x_2)}\sum_{s\in I_q (\Pi_j(t,x),\Pi_j(t,y) )} 
\{ p_{i_1 j}(x_1) p_{i_2 j}(x_2)\overline{p}_s\big ( \Pi_j(t,x_1)\big )\overline{p}_s\big (\Pi_j(t,x_2)\big)\}\, dt \\
&>p_0 q_0 r>0,
\end{aligned}
$$
so {\bf A3} is satisfied. Since $F=X\times X$, assumption {\bf A4} is trivially satisfied.\newline
From Theorem \ref{KS} we obtain (i) and (ii). Set $v(x,i)=C(\varrho (x,x_*)+1)$ and $h(x,i)=C(\varrho(x,x_*)+1+\frac{b}{1-a})$ for $x\in X$, $i\in I$, with $a,b$ as in (\ref{constants}). Iterating (\ref{Lap}) we obtain 
$$
\mathbb{E}_{x_0,_0}(v(x_n,\xi_n))\le h(x_0,\xi_0)\quad\text{for}\quad x_0\in X,\xi_0\in I.
$$
Application of Theorem \ref{Shir} ends the proof.\newline\newline

The next result describing the asymptotic behavior of the process $(x_n)_{n\ge 0}$ on Y is an obvious consequence of Theorem \ref{HS1}. 
Let ${\tilde \mu}_0$ be the distribution of the initial random vector $x_0$ and ${\tilde \mu}_n$ the distribution of $x_n$, i.e.
$$
{\tilde \mu}_n(A)=\mathbb{P}(x_n\in A)\quad\text{for}\quad A\in\mathcal{B}_Y, n\ge 1.
$$
\begin{theorem}
Under the hypotheses of Theorem \ref{HS1} the following statements hold:\newline
(i) there exists a measure ${\tilde \mu}_*\in\mathcal{M}_1^1(Y)$ such that for any ${\tilde \mu}_0$ the sequence 
$({\tilde \mu}_n)_{n\ge 0}$ converges weakly to ${\tilde \mu}_*$. Moreover, if 
$$
\mathbb{P}(x_0\in A)={\tilde\mu}_*(A)\quad\text{for}\quad A\in\mathcal{B}_Y
$$
then ${\tilde\mu}_n(A)={\tilde\mu}_*(A)$ for $A\in\mathcal{B}_Y$ and $n\ge 1$.\newline
(ii)  there exist $q\in (0,1)$ and $C>0$ such that 
$$
||{\tilde\mu}_n - {\tilde\mu}_* ||_{FM}\le q^n C(1+\int_Y \varrho (x,x_*)\,{\tilde\mu}_0(dx))
$$
for any initial distribution ${\tilde\mu}_0\in\mathcal{M}_1^1(Y)$ and $n\ge 1$.\newline
(iii) for any starting point $x_0\in Y$, $\xi_0\in I$ and any bounded Lipschitz function $f$ on $Y$ 
$$
\lim_{n\to\infty}\frac{1}{n}\sum_{k=0}^{n-1} f(x_k)=\int_Y f(x)\,{\tilde\mu}_*(dx)
$$
$\mathbb{P}_{x_0,\xi_0}$ almost surely.
\end{theorem}

\end{document}